\documentclass[11pt]{article}

\usepackage{amsmath,amsthm,amssymb}

\newcommand{\abs}[1]{\left|#1\right|}
\newcommand{\norm}[1]{\left\|#1\right\|}
\newcommand{\normlkp}[3]{\norm{#1}_{L^{#2,#3}}}

\newtheorem{theorem}{Theorem}[section]
\newtheorem{corollary}[theorem]{Corollary}
\newtheorem{lemma}[theorem]{Lemma}
\newtheorem{proposition}[theorem]{Proposition}

\begin{document}

\title{$L^p$ Bernstein Inequalities and Inverse Theorems for RBF Approximation on $\mathbb{R}^d$
\footnote{This paper consists of work from the author's dissertation written under the supervision of Professors F. J. Narcowich and J. D. Ward at Texas A\&M University.  This research was supported by grant DMS-0807033 from the National Science Foundation.
}
}

\author{John Paul Ward}

\date{}

\maketitle

\begin{abstract}
Bernstein inequalities and inverse theorems are a recent development in the theory of radial basis function (RBF) approximation. The purpose of this paper is to extend what is known by deriving $L^p$ Bernstein inequalities for RBF networks on $\mathbb{R}^d$.  These inequalities involve bounding a Bessel-potential norm of an RBF network by its corresponding $L^p$ norm in terms of the separation radius associated with the network.  The Bernstein inequalities will then be used to prove the corresponding inverse theorem.
\end{abstract}

\section{Introduction}
When analyzing an approximation procedure, there are typically two estimates that one is interested in determining.  The first is a direct theorem that gives the rate at which a function may be approximated in terms of its smoothness, and the second is an inverse estimate that guarantees a certain amount of smoothness of an approximant based on its rate of approximation.  Both are equally important, and if the results match up appropriately, they can be combined to completely characterize smoothness spaces in terms of the approximation procedure. In this paper, our goal will be to prove an inverse theorem for RBF approximation.  The usual way one does this is by first deriving a Bernstein inequality, since a standard technique can then be applied to prove the inverse theorem.

Bernstein inequalities date back to 1912 when S.N. Bernstein proved the first inequality of this type for $L^\infty$ norms of trigonometric polynomials, \cite{r:bernstein}.  A generalization can be found in \cite{r:devore2};  this result, which is credited to Zygmund, states that any trigonometric polynomial $T$ of degree $n$ satisfies
\begin{equation*}
 \norm{T^{(r)}}_{L^p} \leq n^r \norm{T}_{L^p}
\end{equation*}
for $1 \leq p \leq \infty$.  However, the first example of a Bernstein-type inequality for RBF approximants was not proved until 2001, \cite{r:schaback}.  Then in 2006, Narcowich, Ward, and Wendland derived a more standard type of Bernstein inequality, \cite{r:narcowich}.  They proved $L^2$ Bernstein inequalities for approximants coming from an RBF approximation space on $\mathbb{R}^d$ where the Fourier transform of the RBF has algebraic decay.  In the same year, Mhaskar proved $L^p$ Bernstein inequalities for certain Gaussian networks on $\mathbb{R}^d$, \cite{r:mhaskar1}, and lastly, in \cite{r:mhaskar2}, Mhaskar, Narcowich, Prestin, and Ward were able to prove Bernstein inequalities in $L^p$ norms for a large class of spherical basis functions.

In this paper we will be concerned with RBF approximants to functions in $L^p(\mathbb{R}^d)$ for $1 \leq p \leq \infty$.  The approximants will be finite linear combinations of translates of an RBF $\Phi$, and the translates will come from a countable set $X \subset \mathbb{R}^d$.  The error of this approximation, which is measured in a Sobolev-type norm, depends on both the function $\Phi$ and the set $X$.  Therefore, given an RBF $\Phi$ and a set $X$, we define the RBF approximation space $S_X(\Phi)$ by 
\begin{equation*}
 S_X(\Phi) = \left\{ \sum_{\xi \in Y} a_\xi \Phi(\cdot-\xi): Y \subset X, \#Y<\infty\right\}\cap L^1(\mathbb{R}^d),
\end{equation*}
and note that these approximation spaces are closely related to the ones studied in 
\cite{r:devore1,r:ward}
, where approximation rates are derived.  By choosing $\Phi$ and $X$ properly, one is able to prove results about rates of approximation as well as the stability of the approximation procedure.  

Our goal will be to establish $L^p$ Bernstein inequalities for certain RBF approximation spaces $S_X(\Phi)$, and these inequalities will take the form $\normlkp{g}{k}{p} \leq Cq_X^{-k}\norm{g}_{L^p}$, where $L^{k,p}$ is a Bessel-potential space.  To prove this, we will use band-limited approximation with the bandwidth proportional to $1/q_X$.  Thus $1/q_X$ acts similarly to a Nyquist frequency, and viewing $1/q_X$ as a frequency, we can see the connection to the classical Bernstein inequalities for trigonometric polynomials.  In particular, bandwidth is playing the role of the degree of the polynomial from the classical inequality.

The RBFs that we will be concerned with have finite smoothness; examples include the Sobolev splines and thin-plate splines. The sets $X$ will be discrete subsets of $\mathbb{R}^d$ with no accumulation points.  For the inverse theorem, we will additionally require that there do not exist arbitrarily large regions with no point from $X$.

\subsection{Strategy}
The basic strategy that we will use is the following, which is the same as the one used in \cite{r:mhaskar2}. Given $g=\sum_{\xi \in X}a_\xi\Phi(\cdot-\xi) \in S_X(\Phi) $, we choose an appropriate band-limited approximant $g_\sigma$, and we have
\begin{equation*}
 \normlkp{g}{k}{p} \leq \normlkp{g_\sigma}{k}{p} + \normlkp{g-g_\sigma}{k}{p}.
\end{equation*}
We then split the second term into two ratios.
\begin{equation}\label{eq:init}
 \normlkp{g}{k}{p} \leq \normlkp{g_\sigma}{k}{p} + \left( \frac{\norm{a}_{\ell_p}}{\norm{g}_{L^p}}\frac{\normlkp{g-g_\sigma}{k}{p}}{\norm{a}_{\ell_p}}\right)\norm{g}_{L^p}.
\end{equation}
The term $\norm{a}_{\ell_p}/\norm{g}_{L^p}$ will be bounded by a stability ratio $\mathcal{R}_{S,p}$ that is independent of the function $g$.  We will then need to bound the error of approximating $g$ by band-limited functions. Combining these results with a Bernstein inequality for band-limited functions, we will be able to prove the Bernstein inequality for all functions in $S_X(\Phi)$, and afterward the corresponding inverse theorem will follow.

\subsection{Notation and Formulas}
For any approximation procedure, one would like to determine the error of the approximation and the stability of the procedure.  When considering an RBF approximation space $S_X(\Phi)$, these quantities are bounded in terms of certain measurements of the set $X$.  The error of approximation is given in terms of the fill distance
\begin{equation*}
   h_X = \sup_{x\in\mathbb{R}^d} \inf_{\xi\in X}\abs{x-\xi},
\end{equation*}
which measures how far a point in $\mathbb{R}^d$ can be from $X$, and the stability of the approximation is determined by the separation radius 
\begin{equation*}
  q_X = \frac{1}{2}\inf_{\genfrac{}{}{0pt}{}{\xi,\xi'\in X}{\xi \neq \xi'}}\abs{\xi-\xi'},
\end{equation*}
which measures how close two points in $X$ may be.  In order to balance the rate of approximation with the stability of the procedure, approximation is restricted to sets $X$ for which $h_X$ is comparable to $q_X$, and sets for which the mesh ratio $\rho_X:= h_X/q_X$ is bounded by a constant will be called quasi-uniform.

Many of the results in this paper will be proved by working in the Fourier domain, and  we will use the following form of the Fourier transform in $\mathbb{R}^d$:

\begin{equation*}
 \widehat{f}(\omega)=\frac{1}{(2\pi)^{d/2}}\int_{\mathbb{R}^d}f(x)e^{-i\omega \cdot x}dx.
\end{equation*}
If $f$ is a radial function, then there is a function $\varphi:(0,\infty)\rightarrow \mathbb{R}$ such that  $f(x)=\varphi(\abs{x})$, and in this case, the Fourier transform of $f$ is given by 

\begin{equation*}
 \widehat{f}(\omega)=\abs{\omega}^{-(d-2)/2}\int_{0}^{\infty}\varphi(t)t^{d/2}J_{(d-2)/2}(\abs{\omega}t) dt,
\end{equation*}
where $J_{(d-2)/2}$ denotes the order $(d-2)/2$ Bessel function of the first kind, cf. \cite[Theorem 5.26]{r:wendland}.

The function spaces that we will be mainly interested in are the Bessel-potential spaces $L^{k,p}(\mathbb{R}^d)$, which coincide with the standard Sobolev spaces $W^{k,p}(\mathbb{R}^d)$ when $k$ is a positive integer and $1<p<\infty$, cf. \cite[Section 5.3]{r:stein1}.   
The Bessel potential spaces are defined by 
\begin{equation*}
 L^{k,p} = \{f:\widehat{f}=(1+\abs{\cdot}^2)^{-k/2}\widehat{g}, g  \in L^p(\mathbb{R}^d)\}
\end{equation*}
for $1 \leq p \leq \infty$, and they are equipped with the norm
\begin{equation*}
 \normlkp{f}{k}{p} = \norm{g}_{L^p}.
\end{equation*}
For the extremal cases $p=1,\infty$, the relationships between the spaces $L^{k,p}$ and $W^{k,p}$ are more complex. For $d=1$ and $k$ even, the spaces are equivalent; however, when $k$ is odd, neither function space is contained in the other for any $d$, cf. \cite[Section 5.6]{r:stein1}. 

\subsection{Radial Basis Functions}
In order to prove the Bernstein inequalities, we will need to require certain properties of the RBFs involved.  As much of the work will be done in the Fourier domain, we state the constraints in terms of the RBFs' Fourier transforms.  Given a radial function $\Phi:\mathbb{R}^d\rightarrow\mathbb{R}$ with (generalized) Fourier transform $\widehat{\Phi}$, let $\phi:(0,\infty)\rightarrow \mathbb{R}$ be the function defined by  $\widehat{\Phi}(\omega)=\phi(\abs{\omega})$.  We will say a function $\Phi$ is admissible of order $\beta$ if there exist constants $C_1,C_2>0$ and $\beta>d$ such that for all $\sigma \geq 1$ and all $l \leq l_d:=\lceil (d+3)/2 \rceil$, the function $F(t):=\phi(t)(1+t^2)^{\beta/2}$ satisfies

\begin{enumerate}
 \item[(i)] $C_1 \leq F(t) \leq C_2$ for $t \geq 1/2$
 \item[(ii)] ${\displaystyle \abs{\left(\frac{d}{dt}\right)^{l} F(t)} \leq C_2 t^{-l}  }$  for $t \geq 1/2$.
\end{enumerate}
Notice, in particular, that the second condition implies 
\begin{equation*}
\abs{\left(\frac{d}{dt}\right)^{l} F(\sigma t)} \leq C \quad \text{for}\quad t\geq 1/2.
\end{equation*}

Two particular classes of admissible functions are the Sobolev splines and the thin-plate splines.  The Sobolev spline $\Phi$ of order $\beta>d$ is given by
\begin{equation*}
 \Phi=C\abs{\cdot}^{(\beta-d)/2}K_{(d-\beta)/2}(\abs{\cdot}),
\end{equation*}
where $K$ is a modified Bessel function of the third kind.  This function possesses the Fourier transform
\begin{equation*}
 \widehat{\Phi}=(1+\abs{\cdot}^2)^{-\beta/2}.
\end{equation*}
With this definition, one can verify that the Sobolev splines of order $\beta$ are the canonical example of admissible functions of order $\beta$ since
\begin{align*}
F(t) &= \phi(t) (1+t^2)^{\beta/2}\\
  &= 1.
\end{align*}
They fit the theory particularly well as they are Green's functions for the Bessel potential differential operators, which we use to measure smoothness.  

For a positive integer $m>d/2$, the thin-plate splines of order $2m$ take the form
\begin{equation*}
\Phi =
\begin{cases}
\abs{\cdot}^{2m-d}, &  d~ \text{odd} \\
\abs{\cdot}^{2m-d}\text{log}\abs{\cdot}, & d~\text{even},
\end{cases}
\end{equation*}
and possess the generalized Fourier transforms
\begin{equation*}
   \widehat{\Phi}=C\abs{\cdot}^{-2m}.
\end{equation*}
Admissibility of order $2m$ can be verified by analyzing the function
\begin{align*}
F(t) &= \phi(t)(1+t^2)^{m} \\
 &= t^{-2m} (1+t^2)^{m}.
\end{align*}
The derivative $F^{(l)}$ can then be bounded as follows:
\begin{align*}
 \abs{F^{(l)}(t)} &= \abs{\sum_{k=0}^l C_k \left(\frac{d}{dt}\right)^{l-k} t^{-2m} \left(\frac{d}{dt}\right)^{k} (1+t^2)^{m}} \\
 &\leq C\sum_{k=0}^l t^{-2m-(l-k)} t^{2m-k} \\
 & \leq Ct^{-l}.
\end{align*}

\section{Stability}\label{sec:stability}
One of the essential results for proving the Bernstein inequalities is a bound of a stability ratio for $S_X(\Phi)$.  We define the $L^p$ stability ratio $\mathcal{R}_{S,p}$ associated with this collection by
\begin{equation*}
 \mathcal{R}_{S,p} = \sup_{S_X(\Phi) \ni g \neq 0} \frac{\norm{a}_{\ell_p}}{\norm{g}_{L^p}},
\end{equation*}
where $g= \sum_{\xi \in X}a_\xi\Phi(\cdot-\xi)$.  The goal of this section is to bound the stability ratio by $Cq_X^{d/p'-\beta}$ for some $C$ independent of $a$ and $X$, where $p'$ is the conjugate exponent to $p$.  For this section, we will assume we are working with a fixed countable set $X\subset \mathbb{R}^d$ with $0< q_X < 1$ and an admissible function $\Phi$ of order $\beta$.

To begin, fix $Y=\{\xi_j\}_{j=1}^N\subset X$ and $g=\sum_{j=1}^N a_j \phi(\cdot-\xi_j)$.  We will derive a bound for $\norm{a}_{\ell_p}/\norm{g}_{L^p}$ and show the bound is independent of $Y$ and $a$.  The strategy for proving this is as follows.  Let $K$ be a smooth function and define $\widehat{K}_\sigma(\omega)=\widehat{K}(\omega/\sigma)$.  We will then consider the convolutions  $K_\sigma*g(x)=\sum_{j=1}^N a_j K_\sigma*\Phi(x-\xi_j)$.  For an appropriate choice of $\sigma$, the interpolation matrix $(A_\sigma)_{i,j}=K_\sigma*\Phi(\xi_i-\xi_j)$ will be invertible, and the norm of its inverse will be bounded.  Then  $a=A_\sigma^{-1}(K_\sigma*g)|_Y$ and $\norm{a}_{\ell_p} \leq \norm{A_\sigma^{-1}}_{\ell_p} \norm{K_\sigma*g|_Y}_{\ell_p}$.  We will then be left with bounding $\norm{K_\sigma*g|_Y}_{\ell_p}$ in terms of $q_X$ and $\norm{g}_{L^p}$.

\subsection{Convolution Kernels}
We now define the class of smooth functions with which we will convolve $g$.
Consider a Schwartz class functions $K:\mathbb{R}^d \rightarrow \mathbb{R}$ that satisfies:
\begin{itemize}
 \item[(i)] There is a $\kappa:[0,\infty)\rightarrow [0,\infty)$ such that $\widehat{K}(\omega)=\kappa(\abs{\omega})$.
 \item[(ii)] $\kappa(r)=0$ for $r\in[0,1]$ and $\kappa$ is non-vanishing on an open set.
\end{itemize}
Given such a $K$, we define the related family $\{K_\sigma\}_{\sigma \geq1}$ by $\widehat{K}_\sigma(\omega)=\widehat{K}(\omega/\sigma)$.  Note that property (ii) requires each function $K_\sigma$ to have a Fourier transform which is $0$ in a neighborhood of the origin, and as $\sigma$ increases, so does this neighborhood.  The convolution $\Phi*K_\sigma$ will retain this property and allow us to obtain diagonal dominance in $A_\sigma$.

Before moving on, we will need to determine certain bounds on the functions $K_\sigma$.  First we need an $L^\infty$ bound.
\begin{align*}
 \abs{K_\sigma(x)}  &\leq C_d \int_{\mathbb{R}^d} \widehat{K}_\sigma(\omega)d\omega \\
                    &\leq C_d\sigma^d \int_1^\infty \kappa(t)t^{d-1}dt,
\end{align*}
so
\begin{equation}\label{eq:ker_bd_0}
 \abs{K_\sigma(x)} \leq C_{d} \sigma^d.
\end{equation}

Next we will need a bound on $K_\sigma$ for $r=\abs{x}>0$.  Writing $K_\sigma$ as a Fourier integral, we see

\begin{align*}
 \abs{K_\sigma(x)}  &= r^{-(d-2)/2}\abs{\int_\sigma^\infty \kappa(t/\sigma)t^{d/2}J_{(d-2)/2}(rt)dt},
\end{align*}
and by a change of variables, we have
\begin{align*}
\abs{K_\sigma(x)} &=  \sigma^{d/2+1} r^{-(d-2)/2}\abs{\int_1^\infty \kappa(t)t^{d/2}J_{(d-2)/2}(\sigma rt)dt}.  
\end{align*}
Since $K$ is a Schwartz class function, $\kappa$ is also smooth. Additionally, $\kappa$ is required to be $0$ in a neighborhood of the origin. Therefore, $\kappa$ satisfies the conditions of Proposition \ref{pr:four1}, and
\begin{equation}\label{eq:ker_bd_1}
 \abs{K_\sigma(x)} \leq C_{d} \frac{\sigma^{d/2+1-l_d}}{r ^{(d-2)/2 +l_d}}.
\end{equation}

Using \eqref{eq:ker_bd_0} and \eqref{eq:ker_bd_1}, we will prove a bound of $L^{p'}$ norms of linear combinations of translates of $K_\sigma$.  

\begin{proposition}\label{prop:MZ_RT}
 Let $T:\mathbb{R}^N\rightarrow L^1(\mathbb{R}^d) \cap L^\infty(\mathbb{R}^d)$ be the linear operator defined by
\begin{equation*}
 T(\gamma) = \sum_{j=1}^N \gamma_j K_\sigma(x-\xi_j).
\end{equation*}
Then
\begin{equation*}
 \norm{T(\gamma)}_{L^{p'}} \leq C_{d} \sigma^{d/p} \left(1 +\frac{1}{(\sigma q_X)^{(d-2)/2 +l_d}} \right)^{1/p}\norm{\gamma}_{\ell_{p'}}.
\end{equation*}
\end{proposition}

\begin{proof}
 After proving the bound in the cases $p=1$ and $p=\infty$, the result will then follow by the Riesz-Thorin theorem (cf.  \cite[Chapter 5]{r:stein2}
).  First, by Proposition \ref{pr:sums} we have,
\begin{align*}
 \sum_{j=1}^{N}\abs{K_\sigma(x-\xi_j)} &\leq  \norm{K_\sigma}_{L^\infty} + \sum_{\abs{x-\xi_j}\geq q} \abs{K_\sigma(x-\xi_j)} \\
  &\leq C_{d} \left( \sigma^d + \frac{\sigma^{d/2+1-l_d}}{q_X ^{(d-2)/2 +l_d}} \right). \\
\end{align*}
Simplifying the previous expression, we obtain
\begin{equation}\label{eq:ker_bd_2}
 \sum_{j=1}^{N}\abs{K_\sigma(x-\xi_j)} \leq C_{d} \sigma^d \left(1 +\frac{1}{(\sigma q_X)^{(d-2)/2 +l_d}} \right).
\end{equation}
Therefore
\begin{align*}
 \norm{T(\gamma)}_{L^\infty} &= \norm{\sum_{j=1}^N \gamma_j K_\sigma(x-\xi_j)}_{L^\infty} \\
  &\leq  \norm{\sum_{j=1}^N  \abs{K_\sigma(x-\xi_j)}}_{L^\infty} \norm{\gamma}_{\ell_\infty} \\
  &\leq  C_{d} \sigma^d \left(1 +\frac{1}{(\sigma q_X)^{(d-2)/2 +l_d}} \right)\norm{\gamma}_{\ell_\infty}
\end{align*}
Now in the case $p=\infty$,
\begin{align*}
 \norm{T(\gamma)}_{L^1} &= \norm{\sum_{j=1}^N \gamma_j K_\sigma(x-\xi_j)}_{L^1} \\
  &\leq \sum_{j=1}^N \abs{\gamma_j} \norm{K_\sigma}_{L^1} \\
  &\leq C \norm{\gamma}_{\ell_1}. 
\end{align*}
\end{proof}

\subsection{Interpolation matrices}
Here, the interpolation matrices $(A_\sigma)_{i,j}=K_\sigma*\Phi(\xi_i-\xi_j)$ will be shown to be invertible by the following lemma.  In fact, the function $K_\sigma*\Phi$ is positive definite, even in the case where $\Phi$ is only conditionally positive definite. Furthermore, the lemma will provide a bound for the $\ell_p$ norm of $A_\sigma^{-1}$.

\begin{lemma}[{\cite[Lemma 5.2]{r:mhaskar2}}]\label{lem:diag_dom}
Given an $n \times n$ matrix $A$, denote its diagonal part by $D$, and let $F=A-D$.
 If $D$ is invertible and $\norm{D^{-1}F}_{\ell_1}<1$, then $A$ is invertible and $\norm{A^{-1}}_{\ell_1}<\norm{D^{-1}}_{\ell_1}(1-\norm{D^{-1}F}_{\ell_1})^{-1}$.
\end{lemma}

The diagonal entries of $A_\sigma$ are equal to $K_\sigma *\Phi(0)$, and the off diagonal absolute column sums are of the form $\sum_{i\neq j}\abs{K_\sigma *\Phi(\xi_i-\xi_j)}$.  In order to apply the lemma, we must bound the former from below and the latter from above.  First,
\begin{align*}
K_\sigma *\Phi(0) &= C_d\int_\sigma^\infty \kappa(t/\sigma)\phi(t)t^{d-1}dt \\
  &= C_d \sigma^{d-\beta} \int_1^\infty \frac{\kappa(t)}{t^{\beta-d+1}} (\sigma t)^{\beta} \phi(\sigma t) dt\\
&\geq C_{\Phi,d} \sigma^{d-\beta} \int_1^\infty \frac{\kappa(t)}{t^{\beta-d+1}} dt.
\end{align*}
The last inequality can be verified by considering the representation

\begin{equation*}
(\sigma t)^{\beta} \phi(\sigma t) = \frac{(\sigma t)^\beta}{(1+(\sigma t)^2)^{\beta/2}} \left(\phi(\sigma t) (1+(\sigma t)^2)^{\beta/2} \right)
\end{equation*}
and applying the definition of admissibility. It now follows that

\begin{equation}\label{eq:conv_bd_0}
 K_\sigma *\Phi(0) \geq  C_{\Phi,d} \sigma^{d-\beta}.
\end{equation}

Next, we need a bound on $\abs{K_\sigma *\Phi(x)}$ for $x\neq0$.  Since $K_\sigma*\Phi$ has a radial Fourier transform in $L_1(\mathbb{R}^d)$, we can write it as a one dimensional integral. Note that in the following integral $r=\abs{x}$.
\begin{align*}
 \abs{K_\sigma *\Phi(x)} &= C_dr^{-(d-2)/2} \abs{ \int_\sigma^\infty \kappa(t/\sigma)\phi(t)t^{d/2}J_{(d-2)/2}(rt)dt} \\
&= C_d \frac{\sigma^{(d+2)/2-\beta}}{r^{(d-2)/2}} \abs{ \int_1^\infty \frac{\kappa(t)}{t^\beta}(\sigma t)^{\beta}\phi(\sigma t) t^{d/2} J_{(d-2)/2}(\sigma rt)dt}. \\
\end{align*}
Admissibility of $\Phi$ together with the decay of $\kappa$ and its derivatives imply that the integrand satisfies the conditions of Proposition \ref{pr:four1}, so

\begin{equation*}
\abs{K_\sigma *\Phi(x)} \leq C_{\Phi,d} \frac{\sigma^{(d+2)/2-\beta} }{r^{(d-2)/2}(\sigma r)^{l_d}}.
\end{equation*}
With this estimate we can bound the off diagonal absolute column sums of $A_\sigma$. Using Proposition \ref{pr:sums} we have
\begin{equation}\label{eq:conv_bd_1}
 \sum_{i\neq j}\abs{K_\sigma*\Phi(\xi_i-\xi_j)} \leq C_{\Phi,d} \frac{\sigma^{d-\beta}}{(\sigma q_X)^{(d-2)/2+{l_d}}}.
\end{equation}
We are now ready to apply the lemma. Define
\begin{center}
 $M = \max\left\{ 1,  \left( \frac{2C_{\Phi,d}^2} {C_{\Phi,d}^1} \right) ^{1/((d-2)/2+{l_d})}\right\},$
\end{center}
 where the constants $C_{\Phi,d}^1$ and $C_{\Phi,d}^2$ are from \eqref{eq:conv_bd_0} and \eqref{eq:conv_bd_1} respectively.  We then define $\sigma_0=M/q_X$, so
\begin{equation*}
 (K_{\sigma_0} *\Phi(0))^{-1} \sum_{i\neq j}\abs{K_{\sigma_0}*\Phi(\xi_i-\xi_j)} \leq \frac{1}{2}.
\end{equation*}
Therefore
\begin{equation*}
 \norm{A_{\sigma_0}^{-1}}_{\ell_1} \leq C_{\Phi,d} \sigma_0^{\beta-d},
\end{equation*}
and in terms of $q_X$, 
\begin{equation}\label{eq:matrix_bd}
 \norm{A_{\sigma_0}^{-1}}_{\ell_1} \leq C_{\Phi,d} q_X^{d-\beta}.
\end{equation}

As $A_{\sigma_0}$ is self-adjoint the same bound holds for $\norm{A_{\sigma_0}^{-1}}_{\ell_\infty}$.  The Riesz-Thorin interpolation theorem can then be applied to get 
\begin{equation}\label{eq:matrix_bd_1}
 \norm{A_{\sigma_0}^{-1}}_{\ell_p} \leq C_{\Phi,d} q_X^{d-\beta}
\end{equation}
for $1 \leq p \leq \infty$.

\subsection{Marcinkiewicz-Zygmund type inequality}
To finish the bound of the stability ratio we require
a bound of a discrete norm by a continuous one. To accomplish this,
we can use an argument similar to the proof of \cite[Theorem 1]{r:mhaskar3}.

\begin{proposition}\label{pr:m-z}
If $1\leq p \leq \infty$ and $f\in L^p$, then
\begin{equation*}
 \norm{K_{\sigma_0}*f|_Y}_{\ell_p} \leq C_{d} q_X^{-d/p} \norm{f}_{L^p}.
\end{equation*}
\end{proposition}
\begin{proof}
Let $p'$ be the conjugate exponent to $p$, i.e. $1/p+1/p'=1$.  
Then since $\ell_{p'}$ is dual to $\ell_{p}$, there is a vector $\gamma\in \mathbb{R}^n$ such that $\norm{\gamma}_{\ell_{p'}}=1$ and 
\begin{equation*}
 \norm{K_{\sigma_0}*f|_{Y}}_{\ell_p}=\sum_{j=1}^{N}\gamma_j K_{\sigma_0}*f(\xi_j).
\end{equation*}
An explicit construction of $\gamma$ can be found in \cite[Proposition 6.13]{r:folland}. 
Writing the convolution as an integral, we get
\begin{align*}
\norm{K_{\sigma_0}*f|_{Y}}_{\ell_p} &= \sum_{j=1}^{N}\gamma_j \int_{\mathbb{R}^d} K_{\sigma_0}(\xi_j-y)f(y)dy \\
& \leq \int_{\mathbb{R}^d} \abs{ \sum_{j=1}^{N}\gamma_j K_{\sigma_0}(\xi_j-y)} \abs{f(y)}dy \\
\end{align*}
We can now apply H\"older's inequality to get 
\begin{equation*}
\norm{K_{\sigma_0}*f|_{Y}}_{\ell_p} \leq  \norm{\sum_{j=1}^{N}\gamma_j K_{\sigma_0}(\cdot-\xi_j)}_{L^{p'}} \norm{f}_{L^p},
\end{equation*}
and finally applying Proposition \ref{prop:MZ_RT} gives the result
\begin{align*}
\norm{K_{\sigma_0}*f|_{Y}}_{\ell_p} &\leq  C_{d} \sigma_0^{d/p} \left(1 +\frac{1}{(\sigma_0 q_X)^{(d-2)/2 +l_d}} \right)^{1/p}\norm{\gamma}_{\ell_{p'}} \norm{f}_{L^p}\\
&\leq C_{d} q_X^{-d/p} \norm{f}_{L^p}.
\end{align*}
\end{proof}

\subsection{Stability Ratio Bound}
We are now in a position to prove the bound on the stability ratio for $p\in[1,\infty]$.  Recall $X$ is a countable subset of $\mathbb{R}^d$ with $0<q_X<1$, and $\Phi$ is an admissible function of order $\beta$.
\begin{theorem}\label{thm:sr}
 Let $\mathcal{R}_{S,p}$ be the stability ratio associated with $S_X(\Phi)$.  Then 
\begin{equation*}
 \mathcal{R}_{S,p}=\sup_{S_X(\Phi) \ni g \neq 0} \frac{\norm{a}_{\ell_p}}{\norm{g}_{L^p}}\leq C_{\Phi,d}q_{X}^{d/p'-\beta}.
\end{equation*}
\end{theorem}
\begin{proof}
 It has been shown that the interpolation matrix $(A_{\sigma_0})_{i,j}=K_{\sigma_0} *\Phi(\xi_i-\xi_j)$ is invertible. Therefore $\norm{a}_{\ell_p} \leq \norm{A_{\sigma_0}^{-1}}_{\ell_p} \norm{K_{\sigma_0}*g|_Y}_{\ell_p}$. Using \eqref{eq:matrix_bd_1}, we get
\begin{equation*}
 \norm{a}_{\ell_p} \leq C_{\Phi,d} q_X^{d-\beta} \norm{K_{\sigma_0}*g|_Y}_{\ell_p}.
\end{equation*}
Finally, applying the M-Z inequality gives the result.
\end{proof}

\section{Band-Limited Approximation}
As in the previous section, $X$ will be a fixed countable set with $0 < q_X <1$, and $\Phi$ will be an admissible function of order $\beta$.  At this point, we are left with bounding the two remaining terms of \eqref{eq:init}.  This will require choosing band-limited functions that approximate the elements of $S_X(\Phi)$ and satisfy the Bernstein inequality as well. In particular, given $g \in S_X(\Phi)$ we need to find a band-limited function $g_{\sigma}$ so that
\begin{equation*}
 \frac{\normlkp{g-g_{\sigma}}{k}{p}}{\norm{a}_{\ell_p}} \leq Cq_{X}^{\beta-k-d/p'}
\end{equation*}
for $1\leq p\leq\infty$. Since most of the work will be done in the Fourier domain, we will impose the condition  $k< \beta - d$.

\subsection{Band-Limited Approximants}
We begin by defining a class of band-limited functions.  A function $g\in S_X(\Phi)$ will be convolved with one of these functions in order to define its band-limited approximant.
Consider a Schwartz class function $K:\mathbb{R}^d \rightarrow \mathbb{R}$ that satisfies the following properties:
\begin{itemize}
 \item[(i)]   There is a non-increasing $\kappa:[0,\infty)\rightarrow [0,\infty)$ such that $\widehat{K}(\omega)=\kappa(\abs{\omega})$
 \item[(ii)] $\kappa(\omega)=1$ for $\omega\leq \frac{1}{2}$, and $\kappa(\omega)=0$ for $\omega \geq1$.
\end{itemize}
Note that this $K$ is different from the one introduced in Section \ref{sec:stability}, and in this section $K$ will be of the form described above.  Given such a $K$, we define the family of functions $\{K_\sigma\}_{\sigma \geq1}$ by $\widehat{K}_\sigma(\omega)=\widehat{K}(\omega/\sigma)$.  Band-limited approximants to $g\in S_X(\Phi)$ are then defined by $g_\sigma = K_\sigma * g$. The first thing we must check is that $g_\sigma$ satisfies the Bernstein inequality.  The following lemma addresses this issue.

\begin{lemma}\label{lem:bl_bi_0}
 Let $f\in L^p(\mathbb{R}^d)$, then
\begin{equation*}
 \normlkp{f*K_\sigma}{m}{p} \leq C_{d}\sigma \normlkp{f*K_\sigma}{m-1}{p} 
\end{equation*}
for $1\leq p\leq\infty$ and any positive integer $m$.
\end{lemma}
\begin{proof}
First, notice that we can write
\begin{equation*}
K_\sigma*f = K_{2\sigma}*(K_\sigma*f),
\end{equation*}
so
\begin{align*}
\left[  (1+\abs{\cdot}^2)^{m/2}\widehat{K}_\sigma \widehat{f}\right]^{\vee} &= \left[(1+\abs{\cdot}^2)^{1/2}\widehat{K}_{2\sigma}  (1+\abs{\cdot}^2)^{(m-1)/2}\widehat{K}_\sigma \widehat{f}\right]^{\vee} \\
& = \left[(1+\abs{\cdot}^2)^{1/2}\widehat{K}_{2\sigma} \right]^{\vee}* \left[(1+\abs{\cdot}^2)^{(m-1)/2}\widehat{K}_\sigma \widehat{f}\right]^{\vee}.
\end{align*}
As $K$ is a Schwartz class function, the first function in the convolution is in $L^1$. Likewise,  $\left[(1+\abs{\cdot}^2)^{(m-1)/2}\widehat{K}_\sigma\right]^{\vee}$ is also in $L^1$. Additionally, $f$ being in $L^p$ implies that the second function of the convolution is in $L^p$. Therefore Young's inequality implies

\begin{align*}
 \normlkp{f*K_\sigma}{m}{p} 
&= \norm{\left[  (1+\abs{\cdot}^2)^{m/2}\widehat{K}_\sigma \widehat{f}\right]^{\vee}}_{L^p} \\
&\leq \norm{\left[(1+\abs{\cdot}^2)^{1/2}\widehat{K}_{2\sigma} \right]^{\vee}}_{L^1} \norm{\left[ (1+\abs{\cdot}^2)^{(m-1)/2}\widehat{K}_\sigma \widehat{f}\right]^{\vee}}_{L^p} \\
&= \norm{\left[(1+\abs{\cdot}^2)^{1/2}\widehat{K}_{2\sigma} \right]^{\vee}}_{L^1} \normlkp{f*K_\sigma}{m-1}{p}.
\end{align*}
Now it is known that there exist finite measures $\nu$ and $\lambda$ such that
\begin{equation*}
(1+\abs{x}^2)^{1/2}=\widehat{\nu}(x)+2\pi\abs{x}\widehat{\lambda}(x),
\end{equation*}
 cf. \cite[Chapter 5]{r:stein1}.   
We therefore have 
\begin{align*}
 \norm{\left[(1+\abs{\cdot}^2)^{1/2}\widehat{K}_{2\sigma} \right]^{\vee}}_{L^1} &= \norm{\left[ \widehat{\nu}\widehat{K}_{2\sigma} +2\pi\abs{\cdot}\widehat{\lambda} \widehat{K}_{2\sigma} \right]^{\vee}}_{L^1} \\
&\leq \norm{\nu *K_{2\sigma}}_{L^1} + 2\pi\norm{\lambda *\left[\abs{\cdot} \widehat{K}_{2\sigma} \right]^{\vee}}_{L^1} \\
&\leq C\left(1+\norm{\left[\abs{\cdot} \widehat{K}_{2\sigma} \right]^{\vee}}_{L^1}\right), \\
\end{align*}
and it remains to prove $\norm{\left[\abs{\cdot}\widehat{K}_{2\sigma} \right]^{\vee}}_{L^1} \leq C\sigma$. First
\begin{align*}
 \norm{ \left[\abs{\cdot}\widehat{K}_{2\sigma} \right]^{\vee}}_{L^\infty} &\leq C_d \int_{\mathbb{R}^d}\abs{\omega} \widehat{K}_{2\sigma}(\omega) d\omega \\
&= C_d \int_0^{2\sigma} t \kappa\left(\frac{t}{2\sigma}\right) t^{d-1}dt \\
&= C_d \sigma^{1+d}\int_0^2  \kappa\left(\frac{t}{2}\right) t^d dt, \\
\end{align*}
and hence
\begin{equation}\label{eq:bd_bla_0}
 \norm{ \left[\abs{\cdot} \widehat{K}_{2\sigma} \right]^{\vee}}_{L^\infty} \leq C_{d} \sigma^{1+d}.
\end{equation}
Now, for $\abs{x}=r>0$ 
\begin{align*}
 \abs{ \left[\abs{\cdot} \widehat{K}_{2\sigma} \right]^{\vee}(x)} &= r^{-(d-2)/2} \abs{ \int_0^{2\sigma} t      \kappa\left(\frac{t}{2\sigma}\right)t^{d/2}J_{(d-2)/2}(rt)dt} \\
&= r^{-(d-2)/2} \sigma^{d/2+2}\abs{ \int_0^2 t \kappa\left(\frac{t}{2}\right)t^{d/2}J_{(d-2)/2}(\sigma rt)dt}, \\
\end{align*}
so by Proposition \ref{pr:four2}
\begin{equation}\label{eq:bd_bla_1}
 \abs{ \left[\abs{\cdot} \widehat{K}_{2\sigma} \right]^{\vee}(r)} \leq C_{d} \frac{r^{-(d-2)/2} \sigma^{d/2+2}}{(\sigma r)^{l_d}}.
\end{equation}
Utilizing inequalities \ref{eq:bd_bla_0} and \ref{eq:bd_bla_1} gives
\begin{align*}
 \norm{ \left[\abs{\cdot} \widehat{K}_{2\sigma} \right]^{\vee}}_{L^1} &\leq C_d \sigma^{-d} \norm{ \left[\abs{\cdot} \widehat{K}_{2\sigma} \right]^{\vee}}_{L^\infty} + \int_{\abs{x}\geq \frac{1}{\sigma}} \abs{\left[\abs{\cdot} \widehat{K}_{2\sigma} \right]^{\vee}(x)} dx \\
&\leq  C_{d} \sigma + C_{d} \int_{\frac{1}{\sigma}}^\infty \frac{r^{-(d-2)/2} \sigma^{d/2+2}}{(\sigma r)^{l_d}} r^{d-1} dr \\
&\leq C_{d}\sigma.
\end{align*}
\end{proof}

\begin{corollary}\label{cor:bl_bi}
Define $\sigma_1$ to be $1/q_X$. Then for all $g\in S_X(\Phi)$,
\begin{equation*}
 \normlkp{g*K_{\sigma_1}}{k}{p} \leq C_{\Phi,d} q_X^{-k}\norm{g}_{L^p}.
\end{equation*}
\end{corollary}

\subsection{Approximation Analysis}
Now that we know the band-limited approximants to the elements of $S_X(\Phi)$ satisfy the Bernstein inequality,  we must bound the error of approximation in $L^{k,p}$. We will begin by bounding the approximation error in $L^{k,1}$ and $L^{k,\infty}$ and then use interpolation to obtain the result for all other values of $p$.  In both extremal cases, this reduces to bounding the error of approximating the RBF by band-limited functions.  For $p=1$ this is straightforward; however, the case $p=\infty$ is more involved.

In order to simplify some expressions, we define the functions
\begin{align*}
 E_{\Phi,k} &:= \abs{\left((1+\abs{\cdot}^2)^{k/2}(\Phi-\Phi*K_{\sigma_1})^\wedge\right)^{\vee}}, \\
   F(t)       &:= \phi(t)(1+t^2)^{\beta/2}.
\end{align*}
If we are to bound the error of approximating $\Phi$ by band-limited functions, we will certainly need a point-wise bound of $E_{\Phi,k}$.  Let us begin with an $L^\infty$ bound. 
\begin{align*}
 E_{\Phi,k}(x) &\leq \int_{\mathbb{R}^d}(1-\widehat{K}_{\sigma_1}(\omega))\widehat{\Phi}(\omega)(1+\abs{\omega}^2)^{k/2}d\omega \\
&\leq C_d {\sigma_1}^{d-\beta+k}\int_{1/2}^\infty \frac{(1-\kappa(t))F({\sigma_1} t)t^{d-1}}{((1/{\sigma_1})^2+t^2)^{(\beta-k)/2}}dt.
\end{align*}
Therefore
\begin{equation}\label{eq:app_bd_0}
  E_{\Phi,k}(x) \leq C_{\beta,d} {\sigma_1}^{d-\beta+k}.
\end{equation}
Next, for $\abs{x}=r>0$
\begin{align*}
 E_{\Phi,k}(x) &= r^{-(d-2)/2} \abs{\int_{{\sigma_1}/2}^\infty  \frac{(1-\kappa(t/{\sigma_1}))F(t)}{(1+t^2)^{(\beta-k)/2}}t^{d/2}J_{(d-2)/2}(rt)dt} \\
&= \frac{ {\sigma_1}^{d/2+1-\beta+k}}{r^{(d-2)/2}} \abs{\int_{1/2}^\infty  \frac{(1-\kappa(t))F({\sigma_1} t)}{((1/{\sigma_1})^2+t^2)^{(\beta-k)/2}}t^{d/2}J_{(d-2)/2}({\sigma_1} rt)dt}. 
\end{align*}
Therefore by Proposition \ref{pr:four1}
\begin{equation}\label{eq:app_bd_1}
  E_{\Phi,k}(x) \leq C_{\Phi,d}  \frac{r^{-(d-2)/2} {\sigma_1}^{d/2+1-\beta+k}}{({\sigma_1} r)^{l_d}} .
\end{equation}
With these results we are now able to bound the error of approximation.  

\begin{theorem}\label{thm:app}
 Let $Y=\{\xi_i\}_{i=1}^N$ be a finite subset of a quasi-uniform set $X$.  Given $g=\sum_{j=1}^Na_j\Phi(\cdot-\xi_j) \in S_X(\Phi)$, we have 
\begin{center}
 $\normlkp{g-g*K_{\sigma_1}}{k}{p} \leq C_{\Phi,d} q_X^{\beta-k-d/p'} \norm{a}_{\ell_p}$
\end{center}
 for $1\leq p \leq \infty$.
\end{theorem}
\begin{proof}
 We will show that this holds when $p=1$ and $p=\infty$, and the result will follow from the Riesz-Thorin interpolation theorem. In this case, the operator being interpolated is $T:\ell_p(Y)\rightarrow L^p(\mathbb{R}^d)$, where
\begin{equation*}
T(a) = ({\rm Id}-\Delta)^{k/2}\left( \sum_{j=1}^N a_j \Phi(\cdot-\xi_j) -K_{\sigma_1} * \left(\sum_{j=1}^N a_j \Phi(\cdot-\xi_j)\right)  \right).
\end{equation*}  
 
Letting $g_{\sigma_1}=g*K_{\sigma_1}$, we have
\begin{align*}
 \normlkp{g-g_{\sigma_1}}{k}{\infty} &= \normlkp{\sum_{j=1}^N a_j \Phi(\cdot-\xi_j) - \sum_{j=1}^N a_j \Phi *K_{\sigma_1}(\cdot-\xi_j)}{k}{\infty} \\
&= \norm{\sum_{j=1}^N a_j\left((1+\abs{\cdot}^2)^{k/2}(\Phi-\Phi*K_{\sigma_1})^\wedge\right)^{\vee} (\cdot-\xi_j)}_{L^\infty}, \\
\end{align*}
and using the notation for $E_{\Phi,k}$ above
\begin{equation*}
 \normlkp{g-g_{\sigma_1}}{k}{\infty} \leq \norm{a}_{\ell_\infty} \norm{\sum_{j=1}^N E_{\Phi,k}(\cdot-\xi_j)}_{L^\infty}
\end{equation*}
Now by \eqref{eq:app_bd_0}, \eqref{eq:app_bd_1}, and Proposition \ref{pr:sums}
\begin{align*}
 \sum_{j=1}^N  E_{\Phi,k}(x-\xi_j) &\leq \norm{E_{\Phi,k} }_{L^\infty} 
 +  \sum_{\abs{x-\xi_j}\geq q_Y}  E_{\Phi,k}(x-\xi_j)   \\
&\leq C_{\Phi,d} {\sigma_1}^{d-\beta+k} \left( 1+ \frac{1}{({\sigma_1} q_X)^{(d-2)/2 +l_d}} \right). \\
\end{align*}
For $p=1$, $\normlkp{g-g_{\sigma_1}}{k}{1} \leq \norm{a}_{\ell_1} \norm{E_{\Phi,k}}_{L^1} $, and 
\begin{align*}
 \norm{E_{\Phi,k}}_{L^1} &\leq C_d q_X^d \norm{E_{\Phi,k}}_{L^\infty} + \int_{\abs{x}\geq q_X} E_{\Phi,k}(x) dx \\
&\leq C_{\Phi,d} {\sigma_1}^{k-\beta} + C_{\Phi,d} \int_q^\infty  \frac{{\sigma_1}^{d/2+1-\beta+k-l_d}}{r^{(d-2)/2 + l_d-d+1}} dr  \\
&\leq C_{\Phi,d} {\sigma_1}^{k-\beta}.
\end{align*}
\end{proof}

\section{Bernstein Inequalities and Inverse Theorems}
In approximation theory, there are a variety of applications for Bernstein inequalities.  While they are most commonly associated with the derivation of inverse theorems, they can also be useful in proving direct theorems. For example, a Bernstein inequality for multivariate polynomials is used in certain RBF approximation error estimates, cf. \cite[Chapter 11]{r:wendland}.    
However, in this paper, we will only address the Bernstein inequalities themselves and their matching inverse theorems.

With the bound of the stability ratio and the band-limited approximation estimate in hand, we are in a position to prove the Bernstein inequalities.  

\begin{theorem}\label{thm:BI}
Let $X$ be a countable set with $0< q_X <1$, and let $\Phi$ be an admissible function of order $\beta$.  If $k<\beta-d$, $1\leq p \leq \infty$, and $g\in S_X(\Phi)$, then 
\begin{equation*}
 \normlkp{g}{k}{p}\leq C_{\Phi,d}q_X^{-k}\norm{g}_{L^p}
\end{equation*}
\end{theorem}
\begin{proof}
Let $g_{\sigma_1}$ be the previously defined approximant of $g$. Then
\begin{align*}
 \normlkp{g}{k}{p} &\leq \normlkp{g_{\sigma_1}}{k}{p} + \normlkp{g-g_{\sigma_1}}{k}{p} \\
&= \normlkp{g_{\sigma_1}}{k}{p} + \left( \frac{\norm{a}_{\ell_p}}{\norm{g}_{L^p}}\frac{\normlkp{g-g_{\sigma_1}}{k}{p}}{\norm{a}_{\ell_p}}\right)\norm{g}_{L^p}
\end{align*}
 Applying Theorem \ref{thm:sr}, Corollary \ref{cor:bl_bi}, and Theorem \ref{thm:app}
\begin{align*}
 \normlkp{g}{k}{p} &\leq \normlkp{g_{\sigma_1}}{k}{p} + \left( \frac{\norm{a}_{\ell_p}}{\norm{g}_{L^p}} \frac{\normlkp{g-g_{\sigma_1}}{k}{p}}{\norm{a}_{\ell_p}}\right)\norm{g}_{L^p} \\
&\leq C_{\Phi,d} q_X^{-k}\norm{g}_{L^p} + \left( C_{\Phi,d}q_X^{\frac{d}{p'} -\beta}\right) \left( C_{\Phi,d} q_X^{\beta-k-\frac{d}{p'}} \right) \norm{g}_{L^p} \\
&\leq C_{\Phi,d}q_X^{-k} \norm{g}_{L^p}.
\end{align*}
\end{proof}

Having established Bernstein inequalities for $S_X(\Phi)$, we can now prove the corresponding inverse theorem. 

\begin{theorem}
Let $\{X_n\}_{n\geq 1}$ be a nested sequence ($X_n \subset X_{n+1}$) of countable sets in $\mathbb{R}^d$ satisfying: $\rho_{X_n}\leq C$ for some constant $C>0$ and $0 < h_{X_n},q_{X_n} < 2^{-n}$.  Furthermore, suppose $1 \leq p \leq \infty$, $f \in L^p(\mathbb{R}^d)$, and $\Phi$ is admissible of order $\beta$. If there is a constant $c_f>0$, independent of $n$, and a positive integer $l$ such that
\begin{equation*}
\inf_{g\in S_{X_n}(\Phi)} \norm{f-g}_{L_p(\mathbb{R}^d)} \leq c_f h_n^l,
\end{equation*}
then $f\in L^{k,p}$ for every $0 \leq k < \min\{\beta-d,l\}$.
\end{theorem}

\begin{proof}
 Let $f_n \in S_{X_n}$ be a sequence of functions satisfying $\norm{f-f_n}_{L^p} \leq 2 c_f h_n^l$. Note that $f_n \in S_{X_m}$ for $m>n$ because the sets $X_n$ are nested. Using the notation $h_n=h_{X_n}$ and $q_n=q_{X_n}$, we have
\begin{align*}
\normlkp{f_{n+1}-f_n}{k}{p}  &\leq  C_{\Phi,d} \left(\frac{h_{n+1}}{q_{n+1}}\right)^k h_{n+1}^{-k}\norm{f_{n+1}-f_n}_{L^p} \\
&\leq  C_{\Phi,d} h_{n+1}^{-k}(\norm{f_{n+1}-f}_{L^p} + \norm{f-f_n}_{L^p}) \\
&\leq  C_{\Phi,d,f} h_{n+1}^{-k}(h_{n+1}^l+h_{n}^l), \\
\end{align*}
and since $h_n<2^{-n}$, it follows that
\begin{equation*}
 \normlkp{f_{n+1}-f_n}{k}{p}  \leq C_{\Phi,d,f} 2^{-(l-k)n}. 
\end{equation*}
This shows $f_n$ is a Cauchy sequence in $L^{k,p}$. Since $L^{k,p}$ is complete, $f_n$ converges to some function $\widetilde{f} \in L^{k,p}$. Since $f_n$ converges to both $f$ and $\widetilde{f}$ in $L^p$, $f=\widetilde{f}$ a.e., and therefore $f \in L^{k,p}$.
\end{proof}

\appendix
\section{Bessel Functions and Fourier Integrals}
\setcounter{theorem}{0}
\renewcommand{\thetheorem}{\Alph{section}.\arabic{theorem}}

A $d$-dimensional Fourier integral of a radial function reduces to a one-dimensional integral involving a Bessel function of the first kind.  Here, we list some of the properties of these Bessel functions and prove bounds for the corresponding Fourier integrals.

\begin{proposition}[{\cite[Proposition  5.4 \& Proposition 5.6]{r:wendland}}]\label{pr:bess} \hfill
 \begin{itemize}
 \item[(1)]  $\frac{d}{dz}\{z^\nu J_\nu(z)\}=z^\nu J_{\nu-1}(z)$
 \item[(2)] $J_{1/2}(z)=\sqrt{\frac{2}{\pi z}} \sin(z)$,  $J_{-1/2}(z)=\sqrt{\frac{2}{\pi z}} \cos(z)$
 \item[(3)] $J_\nu(r)=\sqrt{\frac{2}{\pi r}}\cos(r-\frac{\nu \pi}{2}-\frac{\pi}{4})+\mathcal{O}(r^{-3/2})$ for $r \rightarrow \infty$ and $\nu \in \mathbb{R}$
 \item[(4)] $J_{l/2}^2(r) \leq \frac{2^{l+2}}{\pi r}$ for $r>0$ and $l\in\mathbb{N}$
 \item[(5)] $\lim_{r\rightarrow 0}r^{-l}J_{l/2}^2(r)=\frac{1}{2^l\Gamma^2(l/2+1)}$ for $l\in\mathbb{N}$
\end{itemize}
\end{proposition}

The next proposition makes use of integration by parts in order to bound the Fourier integral of a function whose support lies outside of a neighborhood of the origin.

\begin{proposition}\label{pr:four1}
 Let $\alpha \geq 1$, and let $f\in C^n([0,\infty))$ for some natural number $n>1$. Also, assume there are constants $C,\epsilon>0$ such that $f=0$ on $[0,\frac{1}{2}]$ and $\abs{f^{(j)}(t)}\leq Ct^{-d-\epsilon}$ for $j\leq n$ and $t>1/2$. Then there is a constant $C_\epsilon$ such that \\ $\abs{\int_{1/2}^\infty f(t)t^{d/2}J_{(d-2)/2}(\alpha t)dt} \leq C_\epsilon \alpha^{-n}$
\end{proposition}

\begin{proof}
 We first define a sequence of functions arising when integrating by parts. Let $f_0=f$, $f_1=f'$, and $f_j=\left(\frac{f_{j-1}}{t}\right)'$ for $j\geq2$.  Note that when $j\geq2$, there are constants $c_{j,l}$ such that $f_j(t)=\sum_{l=1}^jc_{j,l}f^{(l)}(t)t^{-2j+l+1}$, and therefore $\abs{f_j(t)}\leq Ct^{-d-j+1-\epsilon}$ for $t>1/2$. Applying the Dominated Convergence Theorem,
\begin{equation*}
\abs{\int_{1/2}^\infty f_0(t)t^{d/2}J_{(d-2)/2}(\alpha t)dt} = \lim_{b\rightarrow \infty} \abs{\int_{1/2}^b f_0(t)t^{d/2}J_{(d-2)/2}(\alpha t)dt}.
\end{equation*}
After integrating by parts and taking the limit we get
\begin{equation*}
 \abs{\int_{1/2}^\infty f_0(t)t^{d/2}J_{(d-2)/2}(\alpha t)dt} = \frac{1}{\alpha} \abs{ \int_{1/2}^\infty \frac{f_1(t)}{t}t^{d/2+1}J_{d/2}(\alpha t)dt}.
\end{equation*}
Integrating by parts $j$ times,
\begin{align*}
 \abs{\int_{1/2}^\infty f_0(t)t^{d/2}J_{(d-2)/2}(\alpha t)dt} &= \frac{1}{\alpha^j}\lim_{b\rightarrow \infty} \abs{\int_{1/2}^b \frac{f_j(t)}{t}t^{d/2+j}J_{d/2+j-1}(\alpha t)dt} \\
&= \frac{1}{\alpha^{j+1}} \abs{\int_{1/2}^\infty \frac{f_{j+1}(t)}{t}t^{d/2+j+1}J_{d/2+j}(\alpha t)dt},
\end{align*}
and therefore
\begin{equation*}
 \abs{\int_{1/2}^\infty f_0(t)t^{d/2}J_{(d-2)/2}(\alpha t)dt} \leq \frac{1}{\alpha^{n}} \int_{1/2}^\infty \abs{\frac{f_{n}(t)}{t}t^{d/2+n}J_{d/2+n-1}(\alpha t)}dt.
\end{equation*}
\end{proof}

The following proposition bounds the Fourier integral of a compactly supported function that is identically one in a neighborhood of the origin.  As in the previous case, the proof relies on integration by parts.

\begin{proposition}\label{pr:four2}
 Let $\alpha \geq 1$, and let $f$ be a function in $C^n([0,\infty))$ for some natural number $n>1$. Also, assume $f=1$ in a neighborhood of $0$ and $f(t)=0$ for $t>2$. Then there is a constant $C$ such that
 \begin{equation*}
 \abs{\int_0^2 tf(t)t^{d/2}J_{(d-2)/2}(\alpha t)dt} \leq C\alpha^{-n}.
 \end{equation*} 
\end{proposition}
\begin{proof}
 We first define a sequence of functions arising when integrating by parts. Let $f_0(t)=tf(t)$, $f_1=f_0'$, and $f_j=\left(\frac{f_{j-1}}{t}\right)'$ for $j\geq2$.  Note that when $j\geq2$, $f_j(t)=O(t^{-2j+2})$ as $t\rightarrow0$. 
After integrating by parts $n$ times, we have
\begin{equation*}
 \abs{\int_0^2 f_0(t)t^{d/2}J_{(d-2)/2}(\alpha t)dt} = \frac{1}{\alpha^{n}} \abs{\int_0^2 \frac{f_{n}(t)}{t}t^{d/2+n}J_{d/2+n-1}(\alpha t)dt}.
\end{equation*}
\end{proof}

\section{Sums of Function Values Over Discrete Sets}
\setcounter{theorem}{0}
\renewcommand{\thetheorem}{\Alph{section}.\arabic{theorem}}

Here we provide a bound for sums of function values taken from discrete sets in $\mathbb{R}^d$. This result is important for showing that the constant in the Bernstein inequality does not depend on the number of centers.

\begin{proposition}\label{pr:sums}
 Let $X \subset \mathbb{R}^d$ be a countable set with $q_X>0$, and let $Y=\{y_j\}_{j=1}^N$ be a subset of $X$ such that $\abs{y_j} \geq q_X$ for $1\leq j\leq N$.  If $f:\mathbb{R}^d \rightarrow \mathbb{R}$ is a function with $\abs{f(x)}\leq C\abs{x}^{-d-\epsilon}$ for some $C,\epsilon>0$, then 
\begin{equation*}
 \sum_{j=1}^N\abs{f(y_j)} \leq 3^d(1+1/\epsilon)(Cq_X^{-d-\epsilon}).
\end{equation*}
\end{proposition}
\begin{proof}
 We can bound the sum using the volume argument found in the proof of \cite[Theorem 12.3]{r:wendland}. 
Following the same procedure, we define 
\begin{equation*}
 E_m=\{x\in\mathbb{R}^d:mq_Y\leq \abs{x} <(m+1)q_Y\}.
\end{equation*}
for each positive integer $m$.  Now by comparing the volume of $E_m$ to the volume of a ball of radius $q_Y$, one finds that $\#Y \cap E_m \leq 3^dm^{d-1}$.  Therefore,
\begin{align*}
 \sum_{j=1}^N \abs{f(y_j)} &\leq \sum_{j=1}^N \frac{C}{\abs{y_j}^{d+\epsilon}}  \\
&\leq \sum_{m=1}^\infty (\#Y \cap E_m) \max_{y\in E_m} \left\{\frac{C}{\abs{y}^{d+\epsilon}}\right\}  \\
&\leq \sum_{m=1}^\infty 3^d m^{d-1} \frac{C}{(mq_Y)^{d+\epsilon}} \\
&\leq 3^d(1+1/\epsilon) Cq_X^{-d-\epsilon}
\end{align*}
\end{proof}

\bibliographystyle{plain}
\bibliography{lp_bernstein}

\begin{thebibliography}{10}

\bibitem{r:bernstein}
S.~N. Bernstein.
\newblock Sur l'ordre de la meilleure approximation des fonctions continues par
  les poly{\^o}mes de degr{\'e} donn{\'e}.
\newblock {\em Mem. Cl. Sci. Acad. Roy. Belg.}, 4:1--103, 1912.

\bibitem{r:devore2}
R.~De{V}ore and G.~G. Lorentz.
\newblock {\em Constructive approximation}, volume 303 of {\em Grundlehren der
  Mathematischen Wissenschaften (Fundamental Principles of Mathematical
  Sciences)}.
\newblock Springer-Verlag, Berlin, 1993.

\bibitem{r:devore1}
R.~De{V}ore and A.~Ron.
\newblock Approximation using scattered shifts of a multivariate function.
\newblock {\em Trans. Amer. Math. Soc.}, 362(12):6205--6229, 2010.

\bibitem{r:folland}
G.~B. Folland.
\newblock {\em Real analysis}.
\newblock Pure and Applied Mathematics (New York). John Wiley \& Sons Inc., New
  York, second edition, 1999.
\newblock Modern techniques and their applications, A Wiley-Interscience
  Publication.

\bibitem{r:mhaskar1}
H.~N. Mhaskar.
\newblock A {M}arkov-{B}ernstein inequality for {G}aussian networks.
\newblock In {\em Trends and applications in constructive approximation},
  volume 151 of {\em Internat. Ser. Numer. Math.}, pages 165--180.
  Birkh{\"a}user, Basel, 2005.

\bibitem{r:mhaskar2}
H.~N. Mhaskar, F.~J. Narcowich, J.~Prestin, and J.~D. Ward.
\newblock {$L^p$} {B}ernstein estimates and approximation by spherical basis
  functions.
\newblock {\em Math. Comp.}, 79(271):1647--1679, 2010.

\bibitem{r:mhaskar3}
H.~N. Mhaskar and J.~Prestin.
\newblock On {M}arcinkiewicz-{Z}ygmund-type inequalities.
\newblock In {\em Approximation theory}, volume 212 of {\em Monogr. Textbooks
  Pure Appl. Math.}, pages 389--403. Dekker, New York, 1998.

\bibitem{r:narcowich}
F.~J. Narcowich, J.~D. Ward, and H.~Wendland.
\newblock Sobolev error estimates and a {B}ernstein inequality for scattered
  data interpolation via radial basis functions.
\newblock {\em Constr. Approx.}, 24(2):175--186, 2006.

\bibitem{r:schaback}
R.~Schaback and H.~Wendland.
\newblock Inverse and saturation theorems for radial basis function
  interpolation.
\newblock {\em Math. Comp.}, 71(238):669--681 (electronic), 2002.

\bibitem{r:stein1}
E.~M. Stein.
\newblock {\em Singular integrals and differentiability properties of
  functions}.
\newblock Princeton Mathematical Series, No. 30. Princeton University Press,
  Princeton, N.J., 1970.

\bibitem{r:stein2}
E.~M. Stein and G.~Weiss.
\newblock {\em Introduction to {F}ourier analysis on {E}uclidean spaces}.
\newblock Princeton University Press, Princeton, N.J., 1971.
\newblock Princeton Mathematical Series, No. 32.

\bibitem{r:ward}
J.~P. Ward.
\newblock ${L}^p$ error estimates for approximation by {S}obolev splines and
  {W}endland functions on $\mathbb{R}^d$.
\newblock {\em Advances in Computational Mathematics}, 2011.
\newblock DOI: 10.1007/s10444-011-9263-7.

\bibitem{r:wendland}
H.~Wendland.
\newblock {\em Scattered data approximation}, volume~17 of {\em Cambridge
  Monographs on Applied and Computational Mathematics}.
\newblock Cambridge University Press, Cambridge, 2005.

\end{thebibliography}

\end{document}